\documentclass [12pt,a4paper] {article}
\usepackage[affil-it]{authblk}

\usepackage [T1] {fontenc}
\usepackage {amssymb}
\usepackage {amsmath}
\usepackage {amsthm}
\usepackage{tikz}
\usepackage{hyperref}
\usepackage{url}

 \usepackage[titletoc,toc,title]{appendix}
 
 \usepackage[margin=1in]{geometry}

\usetikzlibrary{arrows,decorations.pathmorphing,backgrounds,positioning,fit,petri}

\usetikzlibrary{matrix}

\usepackage{enumitem}

\DeclareMathOperator {\td} {td}

\DeclareMathOperator {\ldim} {ldim}
\DeclareMathOperator {\tp} {tp}

\DeclareMathOperator {\ord} {ord}
\DeclareMathOperator {\acl} {acl}

\DeclareMathOperator {\cl} {cl}

\DeclareMathOperator {\DCF} {DCF}

\DeclareMathOperator {\MR} {MR}
\DeclareMathOperator {\MD} {MD}

\DeclareMathOperator {\alg} {alg}

\theoremstyle {definition}
\newtheorem {definition}{Definition} [section]

\newtheorem* {notation} {Notation}

\theoremstyle {plain}

\newtheorem {lemma} [definition] {Lemma}
\newtheorem {theorem} [definition] {Theorem}

\newtheorem {corollary} [definition] {Corollary}

\newtheorem {conjecture} [definition] {Conjecture}

\theoremstyle {remark}
\newtheorem {remark} [definition] {Remark}

\makeatletter
\newcommand {\forksym} {\raise0.2ex\hbox{\ooalign{\hidewidth$\vert$\hidewidth\cr\raise-0.9ex\hbox{$\smile$}}}}
\def\@forksym@#1#2{\mathrel{\mathop{\forksym}\displaylimits_{#2}}}
\def\forkind{\@ifnextchar_{\@forksym@}{\forksym}}

\makeatother

\begin {document}

\title{Ax-Schanuel type theorems and geometry of strongly minimal sets in differentially closed fields}

\author{Vahagn Aslanyan}
\date{}

\maketitle

\begin{abstract}
Let $(K;+,\cdot, ', 0, 1)$ be a differentially closed field. In this paper we explore the connection between Ax-Schanuel type theorems (predimension inequalities) for a differential equation $E(x,y)$ and the geometry of the set $U:=\{ y:E(t,y) \wedge y' \neq 0 \}$ where $t$ is an element with $t'=1$. We show that certain types of predimension inequalities imply strong minimality and geometric triviality of $U$. Moreover, the induced structure on Cartesian powers of $U$ is given by special subvarieties. If $E$ has some special form then all fibres $U_s:=\{ y:E(s,y) \wedge y' \neq 0 \}$ (with $s$ non-constant) have the same properties. 
In particular, since the $j$-function satisfies an Ax-Schanuel theorem of the required form (due to Pila and Tsimerman), our results will give another proof for a theorem of Freitag and Scanlon stating that the differential equation of $j$ defines a strongly minimal set with trivial geometry (which is not $\aleph_0$-categorical though).
\end{abstract}


\section{Introduction}

Throughout the paper we let $\mathcal{K}=(K;+,\cdot, ', 0, 1)$ be a differentially closed field with field of constants $C$. 

Let $E(x,y)$ be (the set of solutions of) a differential equation $f(x,y)=0$ with rational (or, more generally, constant) coefficients. A general question that we are interested in is whether $E$ satisfies an Ax-Schanuel type inequality. A motivating example is the exponential differential equation $y'=yx'$. We know that (the original) Ax-Schanuel (\cite{Ax}) gives a predimension inequality (in the sense of Hrushovski \cite{Hru}) which governs the geometry of our equation. In this case the corresponding reduct of a differentially closed field can be reconstructed by a Hrushovski-style amalgamation-with-predimension construction (\cite{Kirby-semiab}). Zilber calls this kind of predimension inequalities \emph{adequate} (see \cite{Aslanyan-thesis,Aslanyan-adequate-predim} for a precise definition). This means that the reduct satisfies a strong existential closedness property which asserts, roughly speaking, that a system of exponential equations which is not overdetermined has a solution. Being overdetermined can be interpreted as contradicting the Ax-Schanuel inequality. Thus, having an adequate Ax-Schanuel inequality for $E$ will give us a complete understanding of its model theory. For more details on this and related problems see \cite{Aslanyan-thesis,Aslanyan-adequate-predim,Kirby-semiab,Zilb-pseudoexp,Zilb-analytic-pseudoanalytic}. Note however that adequacy will not play any role in this paper, we just mentioned it here in order to give a basic idea why Ax-Schanuel inequalities are important. Ax-Schanuel type statements can also be applied to diophantine geometry. Indeed, they can be used to prove a weak version of the famous Zilber-Pink conjecture in the appropriate setting (see \cite{Zilb-exp-sum-published,Pila-Tsim-Ax-j,Kirby-semiab,Aslanyan-weakMZPD}).

Thus, we want to classify differential equations of two variables with respect to the property of satisfying an Ax-Schanuel type inequality. 
The present work should be seen as a part of that more general project. We explore the connection between Ax-Schanuel type theorems (predimension inequalities) for a differential equation $E(x,y)$ and the geometry of the fibres of $E$. More precisely, given a predimension inequality (not necessarily adequate) for solutions of $E$ of a certain type (which is of the form ``$\td-\dim$'' where $\dim$ is a dimension of trivial type) we show that fibres of $E$ are strongly minimal and geometrically trivial (after removing constant points). Moreover, the induced structure on the Cartesian powers of those fibres is given by \emph{special} subvarieties.

In particular, since an Ax-Schanuel theorem (of the required form) for (the differential equation of) the $j$-function is known (due to Pila and Tsimerman \cite{Pila-Tsim-Ax-j}), our results will give a new proof for a theorem of Freitag and Scanlon \cite{Freitag-Scanlon} stating that the differential equation of $j$ defines a strongly minimal set with trivial geometry (which is not $\aleph_0$-categorical though). In fact, the Pila-Tsimerman inequality is the main motivation for this paper.

Thus we get a necessary condition for $E$ to satisfy an Ax-Schanuel inequality of the given form. This is a step towards the solution of the problem described above. In particular it gives rise to an inverse problem: given a one-variable differential equation which is strongly minimal and geometrically trivial, can we say anything about the Ax-Schanuel properties of its two-variable analogue (see Section \ref{remarks} for more details)?

On the other hand, understanding the structure of strongly minimal sets in a given theory is one of the most important problems in model theory. In $\DCF_0$ there is a nice classfication of strongly minimal sets, namely, they satisfy the Zilber trichotomy (Hrushovski-Sokolovi\'{c} \cite{Hru-Sok}). Hrushovski \cite{Hru-Juan} also gave a full characterisation of strongly minimal sets of order $1$ proving that such a set is either non-orthogonal to the constants or it is trivial and $\aleph_0$-categorical. However there is no general classification of trivial strongly minimal sets of higher order and therefore we do not fully understand the nature of those sets. From this point of view the set $J$ defined by the differential equation of $j$ is quite intriguing since it is the first example of a trivial strongly minimal set in $\DCF_0$ which is not $\aleph_0$-categorical. Before Freitag and Scanlon established those properties of $J$ in \cite{Freitag-Scanlon}, it was mainly believed that trivial strongly minimal sets in $\DCF_0$ must be $\aleph_0$-categorical. The reason for this speculation was Hrushovski's aforementioned theorem on order $1$ strongly minimal sets (and the lack of counterexamples).

Thus, the classification of strongly minimal sets in $\DCF_0$ can be seen as another source of motivation for the work in this paper, where we show that these two problems (Ax-Schanuel type theorems and geometry of strongly minimal sets) are in fact closely related.

Let us give a brief outline of the paper. After defining the appropriate notions we formulate the main results of the paper in Section \ref{j-setup}. We prove those results in Section \ref{j-proofs}. Then we give a brief account of the $j$-function and apply our results to its differential equation in Section \ref{j}. Section \ref{remarks} is devoted to some concluding remarks and inverse questions. We have gathered definitions of several
properties of strongly minimal sets that we need in Appendix \ref{A}. 
\\

\textbf{Acknowledgements.} I would like to thank Boris Zilber and Jonathan Pila for valuable suggestions and comments. 

This work was supported by the University of Oxford Dulverton Scholarship.

\section{Setup and main results}\label{j-setup}
\setcounter{equation}{0}

Recall that $\mathcal{K}=(K;+,\cdot, ', 0, 1)$ is a differentially closed field with field of constants $C$. We may assume (without loss of generality) $\mathcal{K}$ is sufficiently saturated if necessary. Fix an element $t$ with $t'=1$. Let $E(x,y)$ be (the set of solutions of) a differential equation $f(x,y)=0$ with rational (or, more generally, constant) coefficients.

Now we give several definitions and then state the main results of the paper. 



\begin{definition}
Let $\mathcal{P}$ be a non-empty collection of algebraic polynomials $P(X,Y) \in C[X,Y]$. We say two elements $a, b \in K$ are $\mathcal{P}$-\emph{independent} if $P(a,b) \neq 0$ and $P(b,a) \neq 0$ for all $P \in \mathcal{P}$. The $\mathcal{P}$-\emph{orbit} of an element $a \in K$ is the set $\{ b \in K: P(a,b)=0 \mbox{ or } P(b,a)=0 \mbox{ for some } P \in \mathcal{P}\}$ (in analogy with a Hecke orbit, see Section \ref{j}). Also, $\mathcal{P}$ is said to be \emph{trivial} if it consists only of the polynomial $X-Y$.
\end{definition}

\begin{definition}
For a non-constant $ x \in K$ the \emph{differentiation with respect to} $x$ is the derivation $\partial_x : K \rightarrow K $ defined by $\partial_x: y \mapsto \frac{y'}{x'}$.
\end{definition}

Recall that $f(x,y)=0$ is the differential equation defining $E$ and denote $m:= \ord_Yf(X,Y)$ (the order of $f$ with respect to $Y$).

\begin{definition}
We say the differential equation $E(x,y)$ has the $\mathcal{P}$-AS (Ax-Schanuel with respect to $\mathcal{P}$) property if the following condition is satisfied:\\ Let $x_1,\ldots,x_n,y_1,\ldots,y_n$ be non-constant elements of $K$ with $f(x_i,y_i)=0$. If the $y_i$'s are pairwise $\mathcal{P}$-independent then 
\begin{equation}\label{general-Ax-ineq}
\td_CC\left(x_1, y_1, \partial_{x_1}y_1, \ldots, \partial^{m-1}_{x_1}y_1, \ldots, x_n, y_n, \partial_{x_n}y_n, \ldots, \partial^{m-1}_{x_n}y_n\right) \geq mn+1.
\end{equation}
\end{definition}

The $\mathcal{P}$-AS property can be reformulated as follows: for any non-constant solutions $(x_i,y_i)$ of $E$ the transcendence degree in \eqref{general-Ax-ineq} is strictly bigger than $m$ times the number of different $\mathcal{P}$-orbits of $y_i$'s. Note that \eqref{general-Ax-ineq} is motivated by the known examples of Ax-Schanuel inequalities (\cite{Ax,Pila-Tsim-Ax-j,Aslanyan-linear}).

\begin{remark}\label{PvsPredim}
Having the $\mathcal{P}$-AS property for a given equation $E$ may force $\mathcal{P}$ to be ``closed'' in some sense. Firstly, $X-Y$ must be in $\mathcal{P}$. Secondly, if $P_1, P_2 \in \mathcal{P}$ then $P_1(y_1,y_2)=0,~ P_2(y_2,y_3)=0$ impose a relation on $y_1$ and $y_3$ given by $Q(y_1,y_3)=0$ for some polynomial $Q$. Then the $\mathcal{P}$-AS property may fail if $Q\in \mathcal{P}$. 
In that case one has to add $Q$ to $\mathcal{P}$ in order to allow the possibility of an Ax-Shcanuel property with respect to $\mathcal{P}$. 

Similar conditions on $\mathcal{P}$ are required in order for $\mathcal{P}$-independence to define a dimension function (number of distinct $\mathcal{P}$-orbits) of a pregeometry (of trivial type), which would imply that the $\mathcal{P}$-AS property is a predimension inequality. Note that the collection of modular polynomials has all those properties. However, the shape of $\mathcal{P}$ is not important for our results since we assume that a given equation $E$ has the $\mathcal{P}$-AS property.
\end{remark}

\begin{definition}
A $\mathcal{P}$-\emph{special variety} (in $K^n$ for some $n$) is an irreducible (over  $C$) component of a Zariski closed set in $K^n$ defined by a finite collection of equations of the form $P_{ik}(y_i,y_k)=0$ for some $P_{ik} \in \mathcal{P}$. For a differential subfield $L \subseteq K$ a \emph{weakly} $\mathcal{P}$-\emph{special variety over} $L$ is an irreducible (over  $L^{\alg}$) component of a Zariski closed set in $K^n$ defined by a finite collection of equations of the form $P_{ik}(y_i,y_k)=0$ and $y_i=a$ for some $P_{ik} \in \mathcal{P}$ and $a \in L^{\alg}$.
For a definable set $V$, a (weakly) $\mathcal{P}$-\emph{special subvariety} (over $L$) of $V$ is an intersection of $V$ with a (weakly) $\mathcal{P}$-special variety (over $L$).\footnote{A $\mathcal{P}$-special variety $S$ may have a constant coordinate defined by an equation $P(y_i,y_i)=0$ for some $P\in \mathcal{P}$. If no coordinate is constant on $S$ then it is said to be a strongly $\mathcal{P}$-special variety.}
\end{definition}

\begin{remark}
If the polynomials from $\mathcal{P}$ have rational coefficients then $\mathcal{P}$-special varieties are defined over $\mathbb{Q}^{\alg}$. Furthermore, if $E$ satisfies the $\mathcal{P}$-AS property then for the set $U:=\{y: f(t,y)=0 \wedge y' \neq 0 \}$ we have $U \cap C(t)^{\alg} = \emptyset$ and so $\mathcal{P}$-special subvarieties of $U$ over $C(t)$ are merely $\mathcal{P}$-special subvarities.
\end{remark}

\begin{notation} For differential fields $L \subseteq K$ and a subset $A \subseteq K$ the differential subfield of $K$ generated by $L$ and $A$ will be denoted by $L \langle A \rangle$.
\end{notation}

Let $C_0 \subseteq C$ be the subfield of $C$ generated by the coefficients of $f$ and let $K_0=C_0 \langle t \rangle = C_0(t)$ be the (differential) subfield of $K$ generated by $C_0$ and $t$ (clearly $U$ is defined over $K_0$). We fix $K_0$ and work over it (in other words we expand our language with new constant symbols for elements (generators) of $K_0$).

Now we can formulate one of our main results (see Section \ref{A} for definitions of geometric triviality and strict disintegratedness).


\begin{theorem}\label{main-thm}
Assume $E(x,y)$ satisfies the $\mathcal{P}$-AS property for some $\mathcal{P}$. Assume further that the differential polynomial $g(Y):=f(t,Y)$ is absolutely irreducible.
Then 
\begin{itemize}
\item $U:=\{y: f(t,y)=0 \wedge y' \neq 0 \}$ is strongly minimal with trivial geometry.
\item If, in addition, $\mathcal{P}$ is trivial then $U$ is strictly disintegrated and hence it has $\aleph_0$-categorical induced structure.
\item All definable subsets of $U^n$ over a differential field $L \supseteq K_0$ are Boolean combinations of weakly $\mathcal{P}$-special subvarieties over $L$. 
\end{itemize}
\end{theorem}

As the reader may guess and as we will see in the proof, this theorem holds under weaker assumptions on $E$. Namely, it is enough to require that \eqref{general-Ax-ineq} hold for $x_1=\ldots=x_n=t$ (which can be thought of as a weak form of the ``Ax-Lindemann-Weierstrass'' property\footnote{To be more precise, we will say that $E$ has the $\mathcal{P}$-ALW property if the inequality \eqref{general-Ax-ineq} is satisfied under an additional assumption $\td_CC(\bar{x})=1$.}). However, we prefer the given formulation of Theorem \ref{main-thm} since the main object of our interest is the Ax-Schanuel inequality (for $E$).

Further, we deduce from Theorem \ref{main-thm} that if $E$ has some special form, then all fibres $E(s,y)$ for a non-constant $s \in K$ have the above properties (over $C_0\langle s \rangle)$.

\begin{corollary}\label{cor}
Let $E(x,y)$ be defined by $P(x, y, \partial_x y, \ldots, \partial_x^m y)=0$ where $P(X,\bar{Y})$ is an irreducible algebraic polynomials over $C$. Assume $E(x,y)$ satisfies the $\mathcal{P}$-AS property for some $\mathcal{P}$ and let $s \in K$ be a non-constant element. Then 
\begin{itemize}
\item $U_s:=\{y: f(s,y)=0 \wedge y' \neq 0 \}$ is strongly minimal with trivial geometry.
\item If, in addition, $\mathcal{P}$ is trivial then any distinct non-algebraic (over $C_0 \langle s \rangle)$ elements are independent and $U_s$ is $\aleph_0$-categorical.
\item All definable subsets of $U_s^n$ over a differential field $L \supseteq C_0\langle s \rangle$ are Boolean combinations of weakly $\mathcal{P}$-special subvarieties over $L$.
\end{itemize}
\end{corollary}

\begin{remark}
Since $U_s \cap C = \emptyset$, in Theorem \ref{main-thm} and Corollary \ref{cor} the induced structure on $U_s^n$ is actually given by \emph{strongly special} subvarieties (over $L$), which means that we do not allow any equation of the form $y_i=c$ for $c$ a constant. In particular we also need to exclude equations of the form $P(y_i,y_i)=0$ for $P \in \mathcal{P}$.
\end{remark}

We also prove a generalisation of Theorem \ref{main-thm}.

\begin{theorem}\label{main-gen}
Assume $E(x,y)$ satisfies the $\mathcal{P}$-AS property and let $p(Y)\in C(t)[Y] \setminus C,~ q(Y)\in C[Y] \setminus C$ be such that the differential polynomial $f(p(Y),q(Y))$ is absolutely irreducible. Then the set $$U_{p,q}:=\{y: E(p(y),q(y)) \wedge y\notin C \}$$ is strongly minimal and geometrically trivial.
\end{theorem}


As an application of Theorem \ref{main-thm} we obtain a result on the differential equation of the $j$-function which was established by Freitag and Scanlon in \cite{Freitag-Scanlon}. To be more precise, let $F(j,j',j'',j''')=0$ be the algebraic differential equation satisfied by the modular $j$-function (see Section \ref{j}).

\begin{theorem}[\cite{Freitag-Scanlon}]\label{j-trivial}
The set $J \subseteq K$ defined by $F(y,y',y'',y''')=0$ is strongly minimal with trivial geometry. Furthermore, $J$ is not $\aleph_0$-categorical.
\end{theorem}

Strong minimality and geometric triviality of $J$ follow directly from Theorem \ref{main-thm} combined with the Ax-Schanuel theorem for $j$ (see Section \ref{j}). Of course the ``furthermore'' clause does not follow from Theorem \ref{main-thm} but it is not difficult to prove. Theorem \ref{main-thm} also gives a characterisation of the induced structure on the Cartesian powers of $J$. Again, that result can be found in \cite{Freitag-Scanlon}.

The proof of Theorem \ref{j-trivial} by Freitag and Scanlon is based on Pila's modular Ax-Lindemann-Weierstrass with derivatives theorem along with Seidenberg's embedding theorem and Nishioka's result on differential equations satisfied by automorphic functions (\cite{Nishioka-automorphic}). They also make use of some tools of stability theory such as the ``Shelah reflection principle''. However, as one may guess, we cannot use Nishioka's theorem (or some analogue of that) in the proof of \ref{main-thm} since we do not know anything about the analytic properties of the solutions of our differential equation. Thus, we show in particular that Theorem \ref{j-trivial} can be deduced from Pila's result abstractly. The key point that makes this possible is stable embedding, which means that if $\mathcal{M}$ is a model of a stable theory and $X \subseteq M$ is a definable set over some $A \subseteq M$ then every definable (with parameters from $M$) subset of $X^n$ can in fact be defined with parameters from $X\cup A$ (see Section \ref{A}). 


Let us stress once more that the set $J$ is notable for being the first example of a strongly minimal set (definable in $\DCF_0$) with trivial geometry that is not $\aleph_0$-categorical. Indeed the aforementioned result of Hrushovski on strongly minimal sets of order $1$ led people to believe that all geometrically trivial strongly minimal sets must be $\aleph_0$-categorical. Nevertheless, it is not true as the set $J$ illustrates.


\section{Proofs of the main results}\label{j-proofs}





\subsection*{Proof of Theorem \ref{main-thm}}

Taking $x_1=\ldots=x_n=t$ in the $\mathcal{P}$-AS property we get the following weak version of the $\mathcal{P}$-ALW property for $U$ which in fact is enough to prove Theorem \ref{main-thm}.

\begin{lemma}\label{Ax-Lindemann}
$\mathcal{P}$-AS implies that for any pairwise $\mathcal{P}$-independent elements $u_1, \ldots, u_n \in U$ the elements $\bar{u}, \bar{u}', \ldots, \bar{u}^{(m-1)}$ are algebraically independent over $C(t)$ and hence over $K_0$.
\end{lemma}

We show that every definable (possibly with parameters) subset $V$ of $U$ is either finite or co-finite. Since $U$ is defined over $K_0$, by stable embedding there is a finite subset $A=\{a_1, \ldots, a_n\} \subseteq U$ such that $V$ is defined over $K_0\cup A$. It suffices to show that $U$ realises a unique non-algebraic type over $K_0\cup A$, i.e. for any $u_1, u_2 \in U \setminus \acl(K_0 \cup A)$ we have $\tp(u_1/K_0 \cup A)=\tp(u_2/K_0\cup A)$. Let $ u \in U \setminus \acl(K_0\cup A)$. We know that $\acl(K_0\cup A)=(K_0\langle A \rangle) ^{\alg} = (K_0(\bar{a}, \bar{a}', \ldots, \bar{a}^{(m-1)}))^{\alg}$. Since $u \notin (K_0\langle A \rangle)^{\alg}$, $u$ is transcendental over $K_0(A)$ and hence it is $\mathcal{P}$-independent from each $a_i$. We may assume without loss of generality that $a_i$'s are pairwise $\mathcal{P}$-independent (otherwise we could replace $A$ by a maximal pairwise $\mathcal{P}$-independent subset). Applying Lemma \ref{Ax-Lindemann} to $\bar{a},u$, we deduce that $u,u',\ldots, u^{(m-1)}$ are algebraically independent over $K_0\langle A\rangle$. Hence $\tp(u/K_0\cup A)$ is determined uniquely (axiomatised) by the set of formulae
$$\{ g(y)=0 \} \cup \{ h(y) \neq 0: h(Y) \in K_0\langle A \rangle\{ Y\},~ \ord (h) < m \}$$
(Recall that $g$ is absolutely irreducible and hence it is irreducible over any field). In other words $g(Y)$ is the minimal differential polynomial of $u$ over $K_0\langle A \rangle$.

Thus $U$ is strongly minimal. {A similar} argument shows also that if $ A \subseteq U$ is a (finite) subset and $ u \in U \cap \acl(K_0A)$ then there is $a \in A$ such that $u \in \acl(K_0a)$. This proves that $U$ is geometrically trivial.

If $\mathcal{P}$ is trivial then distinct elements of $U$ are independent, hence $U$ is strictly disintegrated.

The last part of Theorem \ref{main-thm} follows from the following lemma.

\begin{lemma}\label{Induced-structure}
Every irreducible (relatively) Kolchin closed (over $C(t)$) subset of $U^n$ is a $\mathcal{P}$-special subvariety of $U^n$.
\end{lemma}

\begin{proof}
Let $V \subseteq U^n$ be an irreducible relatively closed subset (i.e. it is the intersection of $U^n$ with an irreducible Kolchin closed set in $K^n$). Pick a generic point $\bar{v}=(v_1, \ldots, v_n) \in V$ and let $W \subseteq K^n$ be the Zariski closure of $\bar{v}$ over $C$. Let $d:=\dim W$ and assume $v_1, \ldots, v_d$ are algebraically independent over $C$. Then $v_i \in (C(v_1, \ldots, v_d))^{\alg}$ for each $i=d+1, \ldots, n$. By Lemma \ref{Ax-Lindemann} each $v_i$ with $i>d$ must be in a $\mathcal{P}$-relation with some $v_{k_i}$ with $k_i \leq d$. Let $P_i(v_i, v_{k_i})=0$ for $i>d$. The algebraic variety defined by the equations $P_i(y_i, y_{k_i})=0,~ i=d+1, \ldots, n,$ has dimension $d$ and contains $W$. Therefore $W$ is a component of that variety and so it is a $\mathcal{P}$-special variety.

We claim that $W \cap U^n = V$. Since $v_1, \ldots, v_d \in U$ are algebraically independent over $C$, by Lemma \ref{Ax-Lindemann} $\bar{v}, \bar{v}', \ldots, \bar{v}^{(m-1)}$ are algebraically independent over $C(t)$. Moreover, the (differential) type of each $v_i,~ i>d,$ over $v_1, \ldots, v_d$ is determined uniquely by an irreducible algebraic equation. Therefore $\tp(\bar{v}/C(t))$ is axiomatised by formulas stating that $\bar{v}$ is Zariski generic in $W$ and belongs to $U^n$. In other words $\bar{v}$ is Kolchin generic in $W \cap U^n$. Now $V$ and $W\cap U^n$ are both equal to the Kolchin closure of $\bar{v}$ inside $U^n$ and hence they are equal.
\end{proof}

Thus definable (over $C(t)$) subsets of $U^n$ are Boolean combinations of special subvarieties. Now let $L \subseteq K$ be an arbitrary differential subfield over which $U$ is defined. Then definable subsets of $U^n$ over $L$ can be defined with parameters from $\tilde{L} = K_0 \cup (U \cap L^{\alg})$ (see Section \ref{A}). Then Lemma \ref{Induced-structure} implies that irreducible Kolchin closed subsets of $U^n$ defined over $\tilde{L}$ are weakly $\mathcal{P}$-special subvarities of $U^n$ over $L$.


Finally, note that since $U$ does not contain any algebraic elements over $C(t)$, the type of any element $u \in U$ over $C(t)$ is isolated by the formula $f(t,y)=0 \wedge y' \neq 0$.

\subsection*{Proof of Theorem \ref{main-gen}}

We argue as above and show that for a finite set $A=\{a_1,\ldots,a_n\} \subseteq U_{p,q}$ there is a unique non-algebraic type over $K_0\langle A \rangle$ realised in $U_{p,q}$. Here we will use full Ax-Lindemann-Weierstrass. 

If $u \in U_{p,q} \setminus (K_0\langle A \rangle)^{\alg}$ then $q(u)$ is transcendental over $K_0(A)$ and so $q(u)$ is $\mathcal{P}$-independent from each $q(a_i)$.  Moreover, we may assume $\{ q(a_1),\ldots,q(a_n) \} $ is $\mathcal{P}$-independent. Then by the $\mathcal{P}$-AS property 
$$\td_CC\left(p(u), q(u), \ldots, \partial^{m-1}_{p(u)}q(u), p(a_i), q(a_i), \ldots, \partial^{m-1}_{p(a_i)}q(a_i)\right)_{i=1,\ldots,n} \geq m(n+1)+1.$$
But then $$\td_CC\left(t,u,u',\ldots,u^{(m-1)},a_i,a_i',\ldots,a_i^{(m-1)}\right)_{i=1,\ldots,n} \geq m(n+1)+1,$$
and hence $u,u',\ldots,u^{(m-1)}$ are algebraically independent over $K_0\langle A \rangle$. This determines the type $\tp(u/K_0A)$ uniquely as required. It also shows triviality of the geometry.

\subsection*{Proof of Corollary \ref{cor}}

Consider the differentially closed field $\mathcal{K}_s = (K;+,\cdot, \partial_s,0,1)$. 
The given form of the differential equation $E$ implies that $U_s$ is defined over $C_0(s)$ in $\mathcal{K}_s$. However, in general it may not be defined over $C_0(s)$ in $\mathcal{K}$, it is defined over $C_0 \langle s \rangle = C_0(s,s',s'',\ldots)$. Since $\partial_ss=1$, we know by Theorem \ref{main-thm} that $U_s$ is strongly minimal in $\mathcal{K}_s$. On the other hand the derivations $\partial_s$ and $'$ are inter-definable (with parameters) and so a set is definable in $\mathcal{K}$ if and only if it is definable in $\mathcal{K}_s$ (possibly with different parameters). This implies that every definable subset of $U_s$ in $\mathcal{K}$ is either finite or co-finite, hence it is strongly minimal.

Further, Theorem \ref{main-thm} implies that $U_s$ is geometrically trivial over $C_0(s)$ in $\mathcal{K}_s$. By Theorem \ref{geometric-triviality}, $U_s$ is also geometrically trivial over $C_0 \langle s \rangle$ in $\mathcal{K}_s$. On the other hand for any subset $A \subseteq U_s$ the algebraic closure of $C_0 \langle s \rangle \cup A$ is the same in $\mathcal{K}$ and $\mathcal{K}_s$. This implies geometric triviality of $U_s$ in $\mathcal{K}_s$.

The same argument (along with the remark after Theorem \ref{geometric-triviality}) shows that the second and the third parts of Corollary \ref{cor} hold as well.

\section{The modular $j$-function}\label{j}

The modular $j$-invariant satisfies the following order three algebraic differential equation:
\begin{equation}\label{dif-eq-j}
F(y,y',y'',y''')=Sy+R(y)(y')^2=0,
\end{equation}
where $S$ denotes the Schwarzian derivative defined by $Sy = \frac{y'''}{y'} - \frac{3}{2} \left( \frac{y''}{y'} \right) ^2$ and $R(y)=\frac{y^2-1968y+2654208}{2y^2(y-1728)^2}$. Let $J$ be the set defined by \eqref{dif-eq-j}. Note that $F$ is not a polynomial but a rational function. In particular constant elements do not satisfy \eqref{dif-eq-j} for $Sy$ is not defined for a constant $y$. We can multiply our equation through by a common denominator and make it into a polynomial equation 
\begin{equation}\label{j-eq-poly}
F^*(y,y',y'',y''')= q(y)y'y'''-\frac{3}{2}q(y)(y'')^2+p(y)(y')^4 =0,
\end{equation}
where $p$ and $q$ are respectively the numerator and the denominator of $R$. Let $J^*$ be the set defined by \eqref{j-eq-poly}. It is not strongly minimal since $C$ is a definable subset. However, as we will see shortly, $J=J^*\setminus C$ is strongly minimal and $\MR(J^*)=1,~ \MD(J^*)=2$. Thus whenever we speak of the formula $F(y,y',y'',y''')=0$ (which, strictly speaking, is not a formula in the language of differential rings) we mean the formula $F^*(y,y',y'',y''')=0 \wedge y' \neq 0$. 

\begin{definition}
For each $N$ we let $\Phi_N(X,Y)$ be the $N$-th \emph{modular polynomial} (see, for example, \cite{Lang-elliptic}). Denote $\Phi:=\{ \Phi_N(X,Y): N>0\}$. We say two elements are \emph{modularly independent} if they are $\Phi$-independent. For an element $a \in K$ its \emph{Hecke orbit} is its $\Phi$-orbit.
\end{definition}

Let us form the two-variable analogue of equation \eqref{dif-eq-j}:
\begin{equation}\label{Two-var-eq}
f(x,y):=F(y,\partial_xy,\partial^2_xy,\partial^3_xy)=0.
\end{equation}
Now we formulate the (differential version of the) Ax-Schanuel theorem for $j$ established by Pila and Tsimerman in \cite{Pila-Tsim-Ax-j}. 

\begin{theorem}[Ax-Schanuel for $j$]\label{Ax-for-j}
Let $z_i, j_i \in K,~ i=1,\ldots,n$ be non-constant elements such that
\begin{equation*}
f(z_i,j_i)=0.
\end{equation*}
If $j_i$'s are pairwise modularly independent then 
\begin{equation}\label{Ax-ineq}
\td_CC(z_1,j_1,\partial_{z_1}j_1,\partial^2_{z_1}j_1,\ldots,z_n,j_n,\partial_{z_n}j_n,\partial^2_{z_n}j_n) \geq 3n+1.
\end{equation}
\end{theorem}

In other words \eqref{Two-var-eq} has the $\Phi$-AS property. Note that here $K$ does not need to be differentially closed, the result holds for arbitrary differential fields (indeed, the Ax-Schanuel theorem is a universal statement). 

As a consequence of Theorems \ref{main-thm} and \ref{Ax-for-j} we get strong minimality and geometric triviality of $J$ (note that $F^*(Y_0,Y_1,Y_2,Y_3)$ is obviously absolutely irreducible as it depends linearly on $Y_3$).

Lemma \ref{Ax-Lindemann} for $j$ is of course a special case of the Ax-Schanuel theorem for $j$. Nevertheless it can also be deduced from Pila's modular Ax-Lindemann-Weierstrass with derivatives theorem (\cite{Pila-modular}) by employing Seidenberg's embedding theorem. Therefore only Pila's theorem is enough to prove strong minimality and geometric triviality of $J$. Moreover, Corollary \ref{cor} shows that all non-constant fibres of \eqref{Two-var-eq} are strongly minimal and geometrically trivial (after removing constant points) and the induced structure on the Cartesian powers of those fibres is given by (strongly) special subvarieties. Note that it is proven in \cite{Freitag-Scanlon} that the sets $F(y,y',y'',y''')=a$ have the same properties for any $a$.

\begin{remark}
To complete the proof of Theorem \ref{j-trivial}, that is, to show that $J$ is not $\aleph_0$-categorical, one argues as follows (see \cite{Freitag-Scanlon}). The Hecke orbit of an element $ j \in J$ is contained in $J$. Therefore $J$ realises infinitely many types over $j$ and hence is not $\aleph_0$-categorical. 
\end{remark}

\section{Concluding remarks}\label{remarks}
The $\mathcal{P}$-AS property states positivity of a  predimension of the form ``$\td - m \cdot d$'', where $d$ is the number of distinct $\mathcal{P}$-orbits.\footnote{Strictly speaking, we do not know this but can assume it is the case. See Remark \ref{PvsPredim}.} Transcendence degree, being the algebraic dimension, is non-locally modular. On the other hand $d$ is a dimension of trivial type. And it is this fact that is responsible for triviality of the geometry of $U$. 


For the exponential differential equation the predimension is of the form ``$\td - \ldim$'' where $\ldim$ stands for the linear dimension. The linear dimension is locally modular non-trivial and accordingly the strongly minimal set $y'=y$ is not trivial; indeed it is non-orthogonal to $C$. Here strong minimality is obvious as the equation has order one. However one may ask whether strong minimality can be deduced from this type of predimension inequalities in general. The answer is no. For example, consider a linear differential equation with constant coefficients $\partial^2_xy - y=0.$
We showed in \cite{Aslanyan-linear} that an Ax-Schanuel statement holds for it.
However $U$ is the set defined by $y''- y=0 $ and the set $y'=y$ is a definable infinite and co-infinite subset (the differential polynomial $Y''-Y$ is absolutely irreducible though). 

An interesting question is whether there are differential equations with the $\mathcal{P}$-AS property with trivial $\mathcal{P}$. As we showed here, if $E(x,y)$ has such a property then the corresponding $U$ must be strongly minimal and strictly disintegrated. There are quite a few examples of this kind of strongly minimal sets in $\DCF_0$. The two-variable versions of those equations will be natural candidates of equations with the required $\mathcal{P}$-AS property. 

For example, the geometry of sets of the form $y'=f(y)$, where $f$ is a rational function over $C$, is well understood. The nature of the geometry is determined by the partial fraction decomposition of $1/f$. As an example consider the equation
\begin{equation}\label{strict-dis-eq}
y'=\frac{y}{1+y}.
\end{equation}
One can show that it defines a strictly disintegrated strongly minimal set (\cite{Mar-dif}). The two variable analogue of this equation is
\begin{equation}\label{strict-dis-eq-2}
\partial_x y=\frac{y}{1+y}.
\end{equation}
But this is equivalent to the equation $\frac{y'}{y} = (x-y)'$. Denoting $z=x-y$ we get the exponential differential equation $y'=yz'$. It is easy to deduce from this that \eqref{strict-dis-eq-2} does not satisfy the $\mathcal{P}$-AS property with any $\mathcal{P}$ (it satisfies a version of the original exponential Ax-Schanuel inequality though and therefore has a predimension inequality which is of the form $\td - \ldim$). Indeed, the fibre of \eqref{strict-dis-eq-2} above $x=t$ is of trivial type but the section by $x=t+y$ is non-orthogonal to $C$. So according to Theorem \ref{main-gen} the equation \eqref{strict-dis-eq-2} does not satisfy any $\mathcal{P}$-AS property. 
Of course, all the sets $y'=f(y)$ can be treated in the same manner and hence they are not appropriate for our purpose. Thus, one needs to look at the behaviour of all the sets $E(p(y),q(y))$, and if they happen to be trivial strongly minimal sets then one can hope for a $\mathcal{P}$-AS inequality.

The classical Painlev\'{e} equations define strongly minimal and strictly disintegrated sets as well. For example, let us consider the first Painlev\'{e} equation $y''=6y^2+t$. Strong minimality and algebraic independence of solutions of this equation were shown by Nishioka in \cite{Nishioka} (note that strong minimality was discovered earlier by Kolchin (see \cite{Mar-dif})). We consider its two-variable version
\begin{equation}\label{Painleve1}
\partial_x^2y=6y^2+x.
\end{equation}
The goal is to find an Ax-Schanuel inequality for this equation. Note that \eqref{Painleve1} does not satisfy the $\mathcal{P}$-AS property with trivial $\mathcal{P}$. Indeed, if $\zeta$ is a fifth root of unity then the transformation $x \mapsto \zeta^2 x,~ y \mapsto \zeta y$ sends a solution of \eqref{Painleve1} to another solution. If one believes these are the only relations between solutions of the above equation, then one can conjecture the following.

\begin{conjecture}[Ax-Schanuel for the first Painlev\'{e} equation]
If $(x_i,y_i),~ i=1, \ldots, n,$ are solutions to the equation \eqref{Painleve1} and $(x_i/ x_j)^5\neq 1$ for $i\neq j$ then
$$\td(\bar{x},\bar{y},\partial_{\bar{x}}\bar{y}) \geq 2n+1.$$
\end{conjecture}

One could in fact replace $x$'s with $y$'s in the condition $(x_i/ x_j)^5\neq 1$ as those are equivalent. 


Nagloo and Pillay showed in \cite{Nagloo-Pillay} that the other generic Painlev\'{e} equations define strictly disintegrated strongly minimal sets as well. So we can analyse relations between solutions of their two-variable analogues and ask similar questions for them too.

In general, proving that a certain equation has the required transcendence properties may be much more difficult than proving that there are equations with those properties. In this regard we believe that ``generic'' (in a suitable sense) equations must satisfy the $\mathcal{P}$-AS property with trivial $\mathcal{P}$. However we do not go into details here and finish with a final remark.

Zilber constructed ``a theory of a generic function'' where the function has transcendence properties analogous to trivial AS property described here (\cite{Zilb-analytic-pseudoanalytic}). He also conjectured that it has an analytic model, i.e. there is an analytic function that satisfies Zilber's axioms and, in particular, the given transcendence properties. Wilkie constructed \emph{Liouville functions} in \cite{Wilkie-Liouville} and showed that they indeed satisfy the transcendence properties formulated by Zilber. Later Koiran \cite{Koiran} proved that Liouville functions satisfy Zilber's existential closedness statement as well. However, those functions do not satisfy any algebraic differential equation, so we cannot translate the result into a differential algebraic language. If there is a differentially algebraic function with similar properties then it may give rise to a differential equation with the desired properties.

\begin{appendices}
\section{On strong minimality}\label{A}

In this appendix we define some standard properties of strongly
minimal sets that are used throughout the paper. We also prove that geometric triviality of a strongly minimal set does not depend on the set of parameters over which our set is defined. It is of course a well known classical result, but we sketch a proof here since we use the result in the proof of Corollary \ref{cor}. For a detailed account of strongly minimal sets and geometric stability theory in general we refer the reader to \cite{Pillay-geometric}. 

Algebraic closure defines a pregeometry on a strongly minimal set. More precisely, if $X$ is a strongly minimal set in a structure $\mathcal{M}$ defined over $A \subseteq M$ then the operator
$$\cl : Y \mapsto \acl(AY) \cap X, \text{ for } Y \subseteq X,$$
is a pregeometry.

\begin{definition}\label{geometry} 
Let $\mathcal{M}$ be a structure and $X \subseteq M$ be a strongly minimal set defined over a finite set $A \subseteq M$.
\begin{itemize}
\item We say $X$ is \emph{geometrically trivial} (over $A$) if whenever $Y \subseteq X$ and $z \in \acl(AY) \cap X$ then $z \in \acl(Ay)$ for some $y \in Y$. In other words, geometric triviality means that the closure of a set is equal to the union of closures of singletons.
\item $X$ is called \emph{strictly disintegrated} (over $A$) if any distinct elements $x_1, \ldots, x_n \in X$ are independent (over $A$).
\item $X$ is called $\aleph_0$-\emph{categorical} (over $A$) if it realises only finitely many $1$-types over $AY$ for any finite $Y \subseteq X$. This is equivalent to saying that  $\acl(AY) \cap X$ is finite for any finite $Y \subseteq X$.
\end{itemize}
\end{definition}

Note that strict disintegratedness implies $\aleph_0$-categoricity and geometric triviality.

\begin{theorem}\label{geometric-triviality}
Let $\mathcal{M}$ be a model of an $\omega$-stable theory and $X \subseteq M$ be as above. If $X$ is geometrically trivial over $A$ then it is geometrically trivial over any superset $B \supseteq A$.
\end{theorem}
\begin{proof}
By expanding the language with constant symbols for elements of $A$ we can assume that $X$ is $\emptyset$-definable. Also we can assume $B=\{b_1, \ldots, b_n\}$ is finite. Let $z \in \acl(BY)$ for some finite $Y \subseteq X$. By stability $\tp(\bar{b}/X)$ is definable over a finite $C \subseteq X$ and we may assume that $C \subseteq \acl(B) \cap X$.  Therefore $z \in \acl(CY)$. By geometric triviality of $X$ (over $\emptyset$) we have $z \in \acl(c)$ for some $c \in C$ or $ z \in \acl(y)$ for some $  y\in Y$. This shows geometric triviality of $X$ over $B$.
\end{proof}

As we saw in the proof all definable subsets of $X^m$ over $B$ are definable over $\acl(B) \cap X$ (which is the \textit{stable embedding} property). The same argument shows that $\aleph_0$-categoricity does not depend on parameters (see also \cite{Nagloo-Pillay0}). Of course this is not true for strict disintegratedness but a weaker property is preserved. Namely if $X$ is strictly disintegrated over $A$ then any distinct non-algebraic elements over $B$ are independent over $B$.

\end{appendices}


\addcontentsline {toc} {section} {Bibliography}
\bibliographystyle {alpha}
\bibliography {ref}

\end{document}